\documentclass[11pt]{amsart}

\usepackage{graphicx, epstopdf, verbatim, fancyhdr, amscd, amssymb} 
\usepackage[left=1.3in,top=1.1in,right=1.3in,bottom=1.1in,head=.1in]{geometry}

\def\part#1{\smallbreak\noindent{\bf #1\rm.}\vskip.02in} 
\def\idem{\rule[1pt]{40pt}{.7pt}\,,\ }

\newtheorem{theorem}{Theorem}[section]
\newtheorem{claim}[theorem]{Claim}
\newtheorem*{theorem1}{Theorem 1}
\newtheorem*{theorem2}{Theorem 2}
\newtheorem*{theorem3}{Theorem 3}
\newtheorem*{actionlemma}{Action Lemma}

\theoremstyle{definition}
\newtheorem{example}[theorem]{Example}%[section]
\newtheorem{remark}[theorem]{Remark}

\newcommand{\remref}[1]{Remark~\ref{rem:#1}}
\newcommand{\exref}[1]{Example~\ref{ex:#1}}
\newcommand{\secref}[1]{Section~\ref{sec:#1}}

\newcommand{\figref}[1]{Figure~\ref{fig:#1}}

\newcommand{\href}[1]{{#1}}

\newcommand{\foot}[1]{\setcounter{footnote}{1}\footnote{\ #1}}

\def\lto{\longrightarrow}
\def\lhookrightarrow{\ensuremath{\lhook\joinrel\relbar\joinrel\rightarrow}}

\def\ra{\lto}
\def\and{\qquad\text{and}\qquad}

\def\={\ = \ }

\def\bz{\mathbb Z}

\def\br{\mathbb R}
\def\bc{\mathbb C}

\def\bff{\mathbb{F}}
\def\bh{\mathbb H}

\def\bl{\mathbb L}

\def\bu{\mathbb U}
\def\bfa{\bff_2^\alpha}
\def\bx{X}

\def\cpone{\bc P^1}

\def\cpinfty{\bc P^\infty}

\def\hpinfty{\bh P^\infty}

\def\bztwo{\bz_2}
\def\bzfour{\bz_4}

\def\id{\textup{id}}

\def\call{\mathcal L}
\def\calp{\mathcal P}

\def\st{\ | \ }
\def\pt{\textup{pt}}

\def\fr{\bff_1^\alpha}
\def\ufr{\bu_1^\alpha}

\def\bfone{\bff_1(\bx)}
\def\bftwo{\bff_2(\bx)}
\def\bfthree{\bff_3(\bx)}

\def\pin {\pi^n(\bx)} 
\def\pione{\pi^1(\bx)} 
\def\pitwo{\pi^2(\bx)} 
\def\pithree{\pi^3(\bx)} 

\def\hk{H_k(\bx)}
\def\hktwo{H_k(\bx;\bztwo)}
\def\hzero{H_0(\bx)}
\def\hone{H_1(\bx)}
\def\sone{S_1(\bx)}
\def\tone{T_1(\bx)}
\def\htwo{H_2(\bx)}

\def\ttwo{T_2(\bx)}
\def\htwoo{H_2^\circ(\bx)}
\def\hthree{H_3(\bx)}

\def\wtwox{w_2(\bx)}

\def\b{B_1} 
\def\z{Z_1}
\def\c{C_2} 

\def\ext{\textup{Ext}}
\def\hom{\textup{Hom}}
\def\tor{\textup{Tor}}

\def\coker{\textup{coker}}
\def\partialr{\partial_r}
\def\partialz{\partial_z}

\def\dt{\raisebox{-.4ex}{\huge$\cdot$\hskip.005in}}

\def\pt{\textup{pt}}
\def\im{\textup{im}}
\def\rk{\textup{rk}}
\def\stab{\bfone_{\alphat}}

\def\int{{\mathbb I}(\bx)}

\def\we{\tau} % {\widehat e}

\def\conj{\rm conj}

\def\alphat{{\widetilde\alpha}}

\newcommand{\cat}[1]{\mathsf{#1}}
\newcommand{\TOP}{\cat{TOP}}
\newcommand{\CG}{\cat{CG}}
\newcommand{\HG}{\cat{HG}}

\begin{document}
%%%%%%%%%%%%%%%

%%%%%%%%%%%%%%%%%%%%%%%%%%%%%%%%% 
%%%%%%%%%     TITLE and ABSTRACT      %%%%%%%%%%
%%%%%%%%%%%%%%%%%%%%%%%%%%%%%%%%% 

\title{Cohomotopy sets of $4$-manifolds}

\author{Robion Kirby, Paul Melvin and Peter Teichner}
\address{Department of Mathematics, U.C. Berkeley, Berkeley, CA 94720}
\email{kirby@math.berkeley.edu}
\email{teichner@math.berkeley.edu \rm (also at Max-Planck Institut in Bonn, Germany)}
\address{Department of Mathematics, Bryn Mawr College, Bryn Mawr, PA 19010}
\email{pmelvin@brynmawr.edu}

\begin{abstract}
Elementary geometric arguments are used to compute the group of homotopy classes of maps from a $4$-manifold $\bx$ to the $3$-sphere, and to enumerate the homotopy classes of maps from $\bx$ to the $2$-sphere.  The former completes a project initiated by Steenrod in the 1940's, and the latter provides geometric arguments for and extensions of recent homotopy theoretic results of Larry Taylor.  These two results complete the computation of all the cohomotopy sets of closed oriented $4$-manifolds and provide a framework for the study of {\it Morse $2$-functions} on $4$-manifolds, a subject that has garnered considerable recent attention. 
\end{abstract}

\maketitle

\parskip 4pt

\vskip-.1in
\vskip-.1in

%%%%%%%%%%%%%%%%%%%%%%%%%%%%%%%%% 
%%%%%%%%%   INTRODUCTION     %%%%%%%%%%
%%%%%%%%%%%%%%%%%%%%%%%%%%%%%%%%% 

Fix a smooth, closed, connected, oriented $4$-manifold $\bx$.  For each positive integer $n$, consider the {\it cohomotopy set} $\pin = [\bx,S^n]$ of free homotopy classes of maps $\bx\to S^n$.  The Pontryagin-Thom construction gives a refinement of Poincar\'e duality
$$
\begin{CD}
\pin @>{\cong}>> \bff_{4-n}(\bx)\\
@V{h^n}VV @VV{h_{4-n}}V \\ 
H^n(X) @>{\cong}>>  H_{4-n}(X)  
\end{CD}
$$
where $\bff_{k}(\bx)$ is the set of closed $k$-dimensional submanifolds of $X$ with a {\em framing} on their normal bundle, up to normally framed bordism in $X \times [0,1]$.  The maps $h^n$ pull back the cohomology fundamental class of $S^n$ whereas the ``forgetful" maps $h_{k}$ use the normal framing to orient the submanifold and then push forward its homology fundamental class.

The purpose of this note is to compute $\pin$ for all $n$.  In fact the only ``interesting" cases are $\pithree \cong \bfone$ (framed links in $\bx$) and $\pitwo \cong \bftwo$ (framed surfaces in $\bx$), computed in Theorems 1 and  2 below; in all other cases $h^n$ is an isomorphism, by classical arguments recalled in Section 0.   The set  $\pithree$ also has a group structure, inherited from the target  group $S^3$, and $\pitwo$ has an action of this group, induced by the action of $S^3$ on $S^2=S^3/S^1$. We will compute both structures geometrically, the group structure (corresponding to disjoint union in $\bfone$) and the action (corresponding to {\it framed link translation} in $\bftwo$, defined in Section \ref{sec:two}).

Recall that the $4$-manifold $\bx$ is {\em odd} if it contains at least one closed oriented surface of odd self-intersection, and otherwise it is {\em even}.  It is {\em spin} if every surface in $\bx$, orientable or not, has even self-intersection, or equivalently the second Stiefel-Whitney class $w_2(\bx)=0$.

%%%%%%%%%%%%%     THEOREM 1     %%%%%%%%%%%%%%%%%%%
\begin{theorem1} \label{thm:one} 
If $\bx$ is odd, then the forgetful map $h_1:\bfone\to \hone$ is an isomorphism.  
If $\bx$ is even, then there is an extension of abelian groups
\begin{equation}
0 \,\lto\, \bztwo \,\lto\, \bfone \overset{h_1}{\,\lto\,} \hone \,\lto\, 0\,,
\tag{{\bf e}} 
\end{equation}
classified by the unique element of $\ext(H_1(X),\bztwo)$ that maps to $w_2(\bx)$ in the universal coefficient sequence 
$
0 \to \ext(H_1(\bx),\bztwo) %\overset{\delta}
\to H^2(\bx;\bztwo) \to \hom(H_2(X),\bztwo) \to 0\,.
$
In particular {\rm({\bf e})} splits, or equivalently
$
 \bfone \cong H_1(X) \oplus\bztwo,
$
 if and only if $\bx$ is spin. 
\end{theorem1}
%%%%%%%%%%%%%%%%%%%%%%%%%%%%%%%%%%%%%%%%%

\noindent
{\bf Remark.} As stated, Theorem 1 is also true for non-compact $4$-manifolds $\bx$.   The proof given in \secref{one} easily adapts to this case.  In fact, as noted by the referee, the smoothness assumption can be dropped as well since topological $4$-manifolds can be smoothed away from a point (in the closed case).

The set $\pithree$ was first investigated by Steenrod \cite{Steenrod}, building on work of Pontryagin \cite{Pontryagin38} and Eilenberg \cite{Eilenberg}, in the more general context of the study of maps from an arbitrary finite complex to a sphere of codimension one.  Steenrod succeeded in  enumerating the elements of $\pithree$, but did not address the question of its group structure.   Recently, this group structure was analyzed by Larry Taylor \cite{Taylor} from a homotopy theoretic point of view, building on work of Larmore and Thomas \cite{LarmoreThomas}. We discuss in Remark~\ref{rem:Taylor} below how Taylor's Theorem 6.2 can be interpreted in a way that it implies our Theorem 1. 

%%%%%%%%%%%%%     THEOREM 2     %%%%%%%%%%%%%%%%
\begin{theorem2} \label{thm:2}
The set $\pitwo \cong \bftwo$ is a ``sheaf of torsors" over the discrete space $\htwoo$ of classes in $\htwo$ of self-intersection zero.  
% \newline \indent 
More precisely, the forgetful map $h_2:\bftwo \to H_2(\bx)$ has image equal to $\htwoo$.  The fiber $\bfa := h_2^{-1}(\alpha)$ over any $\alpha\in\htwoo$ has a transitive $\bfone$-action by {\rm\bf framed link translation} {\rm (defined in Section \ref{sec:two} below)}, with stabilizer equal to the image of the map $\tau_\alpha:\bfthree\to\bfone$ that records {\rm\bf twice} the intersection with any framed surface representing an element in $\bfa$.  Thus $\bfa$ is an $\fr$-torsor, where $\fr := \coker(\tau_\alpha)$.
\end{theorem2}
%%%%%%%%%%%%%%%%%%%%%%%%%%%%%%%%%%%%%%

\begin{comment}
For any oriented surface $F\subset \bx$ representing such a class, and any framing $\nu$ of its normal bundle, the preimage $ h_2^{-1}[F]$ has a transitive $\bfone$-action by {\itb framed link translation}, with stabilizer equal to %$\stab$ 
twice the image of the {\itb framed intersection map} 
$
i_{F_\nu}: \bfthree \lto\bfone
$
$($see \textup{\secref{two}} for the definitions$)$.  Hence $h_2^{-1}[F]$ is in $($non-canonical$\,)$ one-to-one correspondence with the cokernel of $2i_{F_\nu}$, which in fact is independent of $\nu$ $($see \textup{\remref{coker}}$)$.
% for any normal framing $\nu$ of $F$.  
\end{comment}

This computation of $\pitwo$ provides a framework for the study of {\it Morse $2$-functions} on $4$-manifolds, a subject that has garnered considerable recent attention (see for example \cite{ADK}\cite{GayKirby}\cite{Baykur}\cite{Lekili}\cite{AkbulutKarakurt}\cite{Williams}\cite{GayKirby2}\cite{GayKirby3}\cite{GayKirby4}\cite{Williams2}) and that was the original motivation for our work.\foot{A {\it Morse $2$-function} on $\bx$ is a ``generic" map $\bx\to S^2$.  Every smooth, closed, connected, oriented $4$-manifold admits such a fibration, in fact, one with no ``definite folds" and with connected fibers.  Furthermore, the moves generating homotopic fibrations have been determined.  This sets the stage for a new approach to the study of smooth $4$-manifolds.}   

It should be noted that Theorem 2 can be deduced from Taylor's Theorem 6.6 in \cite{Taylor}, which  applies to a general $4$-complex, although the two proofs are completely different.  Our differential topological approach gives a more geometric perspective, and also provides an answer to a question that was left open in \cite[Remark 6.8]{Taylor} regarding the existence of 4-manifolds of type III$_1$; see the end of \secref{prelim} and \exref{Taylor}.  

In fact, there is a homotopy theoretic version of Theorem 2 that applies to a {\it general} CW-complex.   
It is the special case $(G,T)=(S^3, S^1)$ of the following theorem, which in turn follows from the existence of the fiber bundle
$
G/T \to BT \to BG
$
constructed at the beginning of  \secref{three} (where we also explain the necessary translations).

%%%%%%%%%%%%%     THEOREM 3     %%%%%%%%%%%%%%%%
\begin{theorem3}\label{thm:3}
Let $G$ be a topological group  and $T$ be a closed abelian subgroup of $G$. For any CW-complex $X$ there is an ``exact sequence''
\[
[X,T] \overset{\kappa_u}{\lto} [X,G] \lto [X, G/T] \overset{h}{\lto} [X, BT] \lto [X, BG]
\]
in the following sense: A map $X\to BT$ lifts to a map $X\to G/T$ if and only if it becomes null homotopic when composed with the map $BT\to BG$ induced by the inclusion $T\hookrightarrow G$.  Moreover, the natural action of the group $[X,G]$ on the set $[X,G/T]$ is transitive on the fibers of $h$ and the stabilizer of $u:X\to G/T$ equals the image of the homomorphism 
\[
\kappa_u: [X,T] \lto [X,G]\quad , \quad v \, \longmapsto \, \kappa \circ (u \times v) 
\]
which is induced by the continuous map $\kappa : G/T \times T \to G$ defined by $\kappa(gT,t) = gtg^{-1}$.
\end{theorem3}
%%%%%%%%%%%%%%%%%%%%%%%%%%%%%%%%%%%%%%%

If $G$ is a compact connected Lie group with maximal torus $T$, we invite the reader to check that the degree of $\kappa$ equals the order of the corresponding Weyl group.   This gives a satisfying explanation for the factor of $2$ in Theorem 2, as the order of the Weyl group of the Lie group $S^3$.  We suspect that for non-maximal $T$ the map $\kappa$ is null homotopic; it certainly does factor through a manifold of lower dimension. 

%%%%%%%%%%%%%%%%%%%%%%%%%%%%%%%%% 
%%%%%%%%%%    0. PRELIMINARIES       %%%%%%%%%%
%%%%%%%%%%%%%%%%%%%%%%%%%%%%%%%%% 

\setcounter{section}{-1}
\section{Preliminaries}\label{sec:prelim}

The proofs of Theorems 1--3 are given in the correspondingly numbered sections below.  In this preliminary section, we first recall the classical arguments that  
$$
\pin \ \cong \ H^n(\bx) \quad\text{for $n=1$ or $n\ge 4$}\,,
$$
and that all homology classes in $\bx$ below the top dimension are represented by embedded submanifolds (recall that $\bx$ is assumed smooth) -- a fact that we will often appeal to.   We then set the context for our proofs of Theorems 1 and 2, in which $\pithree$ and $\pitwo$ are computed, introducing a notion of ``twisted" homology classes, and discussing a partition of $4$-manifolds into types I, II, III$_1$ and III$_2$ from properties of their intersection forms.

%%%%%%%%%%%%%%%%%%%%%%%%%%%%%%%%% 
%%%%%%%%%    CLASSICAL COMPUTATONS     %%%%%%%
%%%%%%%%%%%%%%%%%%%%%%%%%%%%%%%%% 

\part{Classical computations} 

The generator $v:S^n \to K(\bz,n)$ is a homotopy equivalence for $n=1$, which implies that $h^1:[Y,S^1] \cong [Y,K(\bz,1)] = H^1(Y)$ for any CW-complex $Y$. If $Y$ is an $n$-complex  then $h^n: [Y,S^n] \cong [Y,K(\bz,n)] = H^n(Y)$ because $v$ is $(n+1)$-connected. Thus for $X^4$ we have $\pin = [X,S^n] \cong H^n(\bx)$ when $n=1$ or $n\ge4$, as stated above.   (Note that
the case $n=4$ was first proved by Hopf, and the cases $n>4$ are immediate from general position.)

A similar discussion applies to the Thom class $u:MSO(n)\to K(\bz,n)$, where $MSO(n)$ is the Thom-space of the universal bundle over the oriented Grassmannian $BSO(n)$. The map $u$ is a homotopy equivalence for $n=1$ or $2$, and in general is $(n+2)$-connected\,\cite{Thom}. Hence $u$ induces an isomorphism $[Y,MSO(n)] \cong H^n(Y)$ for any $(n+1)$-complex $Y$, and so for $X^4$ we have $[X, MSO(n) ] \cong H^n(X)$ for all $n>0$.  Now the {\em oriented} Pontryagin-Thom construction gives a diagram like the first of this paper, with $S^n$ replaced by $MSO(n)$ and normal framings replaced by normal orientations. It follows that  for $k<4$, the group $\hk$ is isomorphic to the set $L_k(\bx)$ of closed oriented $k$-dimensional submanifolds of $X$, up to oriented bordism in $X \times [0,1]$. Note that the orientation of $X$ allows us to translate normal orientations on submanifolds into tangential orientations.

Similarly, using the Thom class $MO(n) \to K(\bztwo,n)$ one can conclude that for $k<4$ but $k\ne2$, $\hktwo$ is isomorphic to the set $L_k(\bx,\bztwo)$ of closed $k$-dimensional (unoriented) submanifolds of $X$, up to bordism in $X \times [0,1]$. For $k=2$, there is an exact sequence
\[
0 \,\lto\, \bz \,\lto\, L_k(\bx,\bztwo) \,\lto\, H_2(X;\bztwo) \,\lto\, 0
\]
(see \cite{Thom} and also \cite{FreedmanKirby}) so classes in $H_2(X;\bztwo)$ are still represented by submanifolds. In the following, we will frequently use these interpretations of homology classes in terms of submanifolds in $X$. In particular, the forgetful maps $h_n:\bff_n(X) \to H_n(X)$  are then just given by only remembering the orientation from a framing.  The reduction mod 2 map $H_n(X) \to H_n(X;\bztwo)$ forgets the orientation on the submanifold.

%As noted above, the Pontryagin-Thom construction gives a one-to-one correspondence between the cohomotopy set $\pin$ and the set $\bff_{4-n}(\bx)$, where
%$$
%\bfk \ = \ 
%\{\text{framed bordism classes of framed $k$-dimensional submanifolds of $\bx$}\}.
%$$
%All ``framings" of submanifolds are understood to be trivializations of their normal bundles, which, together with the orientation of $X$, serve to orient these submanifolds.  This correspondence associates to each homotopy class in $\pin$, represented by a smooth map $f:\bx\to S^n$, the framed bordism class of $f^{-1}(y)$ for any regular value $y$ of $f$; see for example Milnor \cite{Milnor65}.  The framed bordisms are submanifolds of $X \times I$ with normal framings that restrict to the given framings on the boundary.    We will refer to $\bfk$ as the {\it $k^\textup{th}$ framed bordism set} of $\bx$.   Our task is to compute $\bfone$ and $\bftwo$
%
%\begin{comment}
%$$
%\begin{aligned}
%\bfone \ &= \ \{\text{framed bordism classes of framed links in $\bx$}\} % \ \text{and}
%\\
%\bftwo \ &= \ \{\text{framed bordism classes of framed surfaces in $\bx$}\}
%\end{aligned}
%$$
%\end{comment}
%
%These computations will be expressed in terms of the forgetful maps
%$$
%h_k : \bfk \lto \hk \quad,\quad [K_\nu]\,\lmapsto\,[K]
%$$
%for $k=1$ and $2$.  Here $K$ is a submanifold (a {\it link} for $k=1$ or a {\it surface} for $k=2$) with normal framing $\nu$, and $h_k$ maps the {\it framed bordism class} of $K_\nu$ to the {\it homology class} (or equivalently {\it bordism class}) of $K$.

\part{Twisted classes}\label{twisted}  

Our study of the framed bordism sets $\bfone$ and $\bftwo$ will feature two special subsets  $\tone\subset\tor_2(\hone)$ and $\ttwo\subset\htwoo$, whose elements we call {\it twisted classes}.  Here $\tor_2(\hone)$ denotes the subgroup of $\hone$ of all elements of order at most $2$, and $\htwoo$ denotes the subset of $\htwo$ of all classes of self-intersection zero. 

To define these subsets, first observe that $\bfone$ is an abelian group with respect to the operation of disjoint union, and that the forgetful map $h_1:\bfone\to\hone$ is an epimorphism.  Each $\kappa\in\hone$ is represented by a knot whose normal bundle has exactly two trivializations (up to homotopy) since $\pi_1SO(3)=\bztwo$.  If the resulting pair of framed knots are framed bordant, then they represent the unique element in the preimage $h_1^{-1}(\kappa)$, and if not, then they represent the two distinct elements in $h_1^{-1}(\kappa)$.  Thus $h_1$ is either an isomorphism or a two-to-one epimorphism (Theorem 1 refines this statement).

Now consider the image subgroup $\sone$ under $h_1$ of the $2$-torsion subgroup of $\bfone$, and define $\tone$ to be its complement in $\tor_2(\hone)$:
$$
\sone \ = \ h_1(\tor_2(\bfone)) \and \tone \ = \ \tor_2(\hone) - \sone.
$$
Evidently $h_1$ splits over $\sone$ (that is, there is a homomorphism $s:\sone\to\bfone$ with $h_1\circ s = \id$) and so we refer to the elements in $\sone$ as {\it split classes}.  The elements in $\tone$, which we call {\it twisted classes}, are exactly those $2$-torsion classes that generate subgroups over which $h_1$ does {\it not} split.  In geometric terms, if a twisted class is represented by a knot $K$, then $K$, equipped with either framing, represents an element of order $4$ in $\bfone$.  Clearly $\tone$ is nonempty if and only if $h_1$ does not split over $\hone$, and in this case it is the nontrivial coset of the index two subgroup $\sone\subset\tor_2(\hone)$.

A class in $\htwoo$ will be called {\it twisted} if its intersection with at least one $3$-dimensional homology class in $\bx$ is a twisted $1$-dimensional class:
$$
\ttwo \ = \ \{\alpha\in\htwoo \st \alpha\dt\tau \in \tone \text{ for some $\tau\in\hthree$}\}.
$$
Evidently $\ttwo$ is empty when $\tone$ is empty, but as will be seen below, the converse may fail.  

%%%%%%%%%%%%%%%%%%%%%%%%%%%%%%%%% 
%%%%%%%%%    4-MANIFOLD TYPES     %%%%%%%%%%%
%%%%%%%%%%%%%%%%%%%%%%%%%%%%%%%%% 

\part{4-manifold types I, II and III}  

The {\it parity} of the  $4$-manifold $\bx$ is determined by the self-intersections of the {\em orientable} surfaces in $\bx$.  If they are all even, then $X$ is said to be {\it even}, and otherwise it is {\it odd}.  The even ones includes all the {\it spin} $4$-manifolds -- those whose tangential structure groups lift to Spin(4) -- which are characterized by the condition that all surfaces in $\bx$, orientable or not, have even self-intersections (only defined modulo $2$ for non-orientable surfaces); by the Wu formula, this is equivalent to the condition $w_2(X)=0$.  As is customary, we say $\bx$ is of {\it type} I, II or III, according to whether it is odd, spin, or even but not spin.   

\begin{figure}[h!]
\includegraphics[height=80pt]{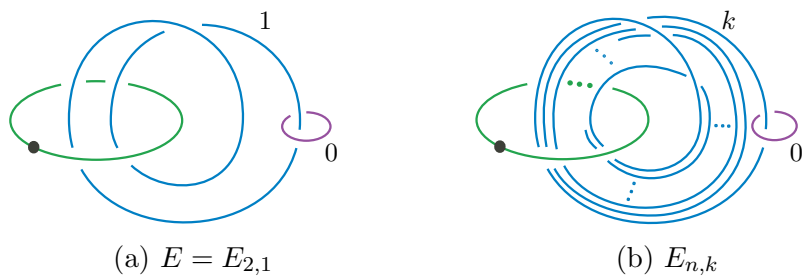}
\put(-260,-15){(a) $E = E_{2,1}$}
\put(-205,75){\small$1$}
\put(-180,25){\small$0$}
\put(-70,-15){(b) $E_{n,k}$}
\put(-30,75){\small$k$}
\put(-4,25){\small$0$}
\caption{The $4$-manifolds $E_{n,k}$}
\label{fig:even}
\end{figure}

All $4$-manifolds of type III must have $2$-torsion in their first homology.  In fact by Theorem 1 (proved in the next section) $\bx$ is of type III if and only if $\hone$ contains twisted classes, i.e.\ $\tone$ is nonempty.  The simplest such manifold is the unique nonspin, oriented $4$-manifold $E$ that fibers over $\br P^2$ with fiber $S^2$.  This manifold is {\it even} since the fiber, which generates $H_2(E;\bz)=\bztwo$, has zero self-intersection, whereas any section has odd self-intersection.  A handlebody picture of $E$, minus the $3$ and $4$-handle, is shown in \figref{even}(a) using the conventions of \cite[Ch.1]{Kirby}.  The generator of $H_1(E)$, represented by the core of the obvious Mobius band bounded by the (attaching circle of the) 1-framed 2-handle, is a twisted $1$-dimensional class.  

%%% FIGURE 1: Non-spin even 4-manifolds %%%
%%%%%%%%%%%%%%%%%% 

Note that $E$ is one of a doubly indexed family of $4$-manifolds $E_{n,k}$ shown in \figref{even}(b), with $H_2(E_{n,k})=\bz_{n}$.  Here $n$ is the linking number between the dotted circle (the $1$-handle) and the $k$-framed 2-handle.   Alternatively, these manifolds can described as the boundaries of $S^1\times B^4$ with a single $2$-handle added.  An easy exercise in link calculus (or even easier $5$-dimensional argument) shows that for fixed $n$, the diffeomorphism type of $E_{n,k}$ depends only on the parity of $k$ when $n$ is even, and is independent of $k$ when $n$ is odd.   Thus there are really only two such manifolds $E_{n,0}$ and $E_{n,1}$ when $n$ is even, and only one $E_{n,0}$ when $n$ is odd.  The $E_{n,0}$'s are all of type II, while the $E_{n,1}$'s (for even $n$) are all of type III.   These manifolds will feature as building blocks for our examples below. 

%%%%%%%%%%%%%%%%%%%%%%%%%%%%%%%%% 
%%%%%%%%%    TYPES III_1 AND III_2     %%%%%%%%%%%
%%%%%%%%%%%%%%%%%%%%%%%%%%%%%%%%% 

\part{Types III$_\text{1}$ and III$_\text{2}$}\label{types}  

When computing their second framed bordism sets, the 4-manifolds $\bx$ of type III will be seen to display the most interesting behavior, suggesting a partition into two subclasses:  $\bx$ is of {\it type} III$_1$ if it contains some twisted 2-dimensional classes, i.e.\ $\ttwo$ is nonempty, and is otherwise of {\it type} III$_2$.  
Thus the 4-manifolds of type III$_1$ are exactly those for which {\it some} of their twisted 1-dimensional classes (which they have because they are of type III) arise as intersections of homology classes of dimensions 2 and 3.    

There exist $4$-manifolds of both subtypes.  In particular, the $E_{n,1}$'s for even $n$ are all of type III$_2$ since they have vanishing third homology.  In \exref{Taylor}, we will construct a $4$-manifold of type III$_1$, which answers a question raised in \cite[Remark 6.8]{Taylor}. It will be shown in \remref{coker} that with regard to their cohomotopy theory, the manifolds of type III$_1$ behave more like odd manifolds, while those of the type III$_2$ behave more like spin manifolds, which explains the choice of subscripts.

%%%%%%%%%%%%%%
\begin{comment}
there is an embedded projective plane in $\bx_1$ (the union of the core of the $1$-framed $2$-handle with the obvious M\"obius band bounded by its attaching circle) of odd self-intersection, and so $\bx_1$ is not spin.
\end{comment}
%%%%%%%%%%%%%%

\medskip
%%%%%%%%%%%%%%%%%%%%%%%%%%%%%%%%% 
%%%    SEC 1: THE FIRST FRAMED BORDISM GROUP    %%%%
%%%%%%%%%%%%%%%%%%%%%%%%%%%%%%%%% 
\section{Computation of the group $\bfone\ \cong \ \pithree$} \label{sec:one}
%%%%%%%%%%%%%%%%%%%%%%%%%%%%%%%%%
%%%%%%%%%%%%%%%%%%%%%%%%%%%%%%%%%

As noted in the introduction, the set $\bfone$ of framed bordism classes of framed links in $\bx$ is an abelian group under the operation of disjoint union.   This operation corresponds to the product in $\pithree$ inherited from the group structure on the target $3$-sphere.  Indeed, given two maps $f_1,f_2:\bx\to S^3$, one can pull back regular values to get disjoint framed links $\bl_1,\bl_2$ in $\bx$.  Choose disjoint tubular neighborhoods $T_1,T_2$ of $\bl_1,\bl_2$.  Up to homotopy, $f_i$ is given by mapping $\bl_i$ to $1$, wrapping the disk fibers in $T_i$ around the $3$-sphere using the framing, and mapping $\bx-T_i$ to $-1$.  It follows that the product $f_1f_2$ has $-1$ as a regular value, with pre-image equal to the (disjoint) union $\bl_1\cup \bl_2$.  This yields an easy indirect proof that $\pithree$ is abelian; there is also a direct proof, well known to homotopy theorists, following from the observation that $S^3$ can be replaced with the fiber of the map $Sq^2:K(\bz,3) \to K(\bztwo,5)$, which is a homotopy-abelian $H$-space (cf.\ \cite[\S6.1]{Taylor}).

It is evident that the forgetful map  
$$
h_1 : \bfone \longrightarrow \hone\,, 
$$
sending the bordism class of a framed link to the homology class of the underlying link, is a surjective homomorphism. 
%%%%%%%%%
\begin{comment}
From the dual perspective, this is just the map 
$$
h^1:\pithree\longrightarrow H^3(\bx;\bz) \quad,\quad [f]\longmapsto f^*(\mu)
$$
where $\mu$ is the orientation class in $H^3(S^3;\bz)$.  Steenrod, in his 1947 paper \cite{Steenrod} introducing his squaring operations, identified the kernel of this map with the quotient group $H^4(\bx,\bz)/Sq^2H^2(\bx;\bz)$.    
\end{comment} 
%%%%%%%%%

\begin{claim}\label{prop:claim} 
The kernel of $h_1$ is either trivial or cyclic of order two, according to whether $\bx$ is an odd or even $4$-manifold.
\end{claim}

This fact was known to Steenrod.  In his 1947 paper \cite{Steenrod} (where he introduced his squaring operations) he identified the kernel of the map $h^1:\pithree\to H^3(\bx)$ dual to $h_1$ with the cokernel of the squaring map $\overline{Sq}^2:H^2(\bx) \to H^4(\bx;\bztwo)$.  Since $\overline{Sq}^2$ is dual to mod 2 self-intersection $\htwo \to H_0(\bx;\bztwo) = \bztwo$ of integral classes, this cokernel is $0$ or $\bztwo$ according to whether $\bx$ is odd or even.  

Thus for odd $\bx$ we have $\pithree\cong H^3(\bx)$, or dually  $\bfone\cong\hone$.  For even $\bx$ there is an abelian extension
$$
0 \lto \bztwo \lto \pithree \overset{h^1}\lto H^3(\bx) \lto 0
$$
\vskip.05in\noindent  (cf.\  \cite[\S6.1]{Taylor}), or dually
\begin{equation}
0 \,\lto\, \bztwo \overset{u_1}{\,\lto\,} \bfone \overset{h_1}{\,\lto\,} \hone \,\lto\, 0
\tag{{\bf e}} 
\end{equation}
as asserted in Theorem 1.  

For completeness, we provide a geometric proof of Claim \ref{prop:claim}, in terms of framed links, which also yields a simple description for the map $u_1$ in $({\bf e})$\,:  Any element in $\ker(h_1)$ is represented by an unknot $U\subset \bx$ with one of its two possible framings.  Exactly one of these framings extends over any given proper $2$-disk $D$ in $\bx\times I$ bounded by $U$.  Let $U_0$ denote the resulting {\it $0$-framed unknot} (this may depend on the choice of $D$) and $U_1$ denote $U$ with the other framing.  Then clearly $[U_0]=0$ in $\bfone$, while $[U_1]$ generates $\ker(h_1)$ and is of order at most two.   The claim is that $[U_0]=[U_1]$ if and only if $X$ is odd.  But if $[U_0]=[U_1]$, then capping off the boundary of a framed bordism in $\bx\times I$ from $U_0$ to $U_1$ and  projecting to $\bx$ produces a (singular) surface of odd self-intersection, so $\bx$ is odd. Conversely, if $\bx$ is odd, then removing two disks from an embedded oriented surface in $\bx$ of odd self-intersection and then ``tilting" this surface in $\bx\times I$ gives a framed bordism between $U_0$ and $U_1$.  This proves the claim, and shows that in the even case, the map $u_1$ in $({\bf e})$ sends the generator of $\bztwo$ to $[U_1]$.  

%%%%%%%%%%
\begin{comment}
Taking $U$ to be the boundary of an embedded $2$-disk $D\subset X$, exactly one of these framings extends over $D$ (pushed into $X\times I$) and the resulting {\it $0$-framed unknot} $U_0$ represents the zero element in the group $\bfone$.  The other framing gives the {\it $1$-framed unknot} $U_1$, which may or may not be framed bordant to $U_0$.   These are pictured  below.
The picture is in an $\br^3$ slice of a coordinate chart, where the tip of the first vector $e_1$ in the framing is drawn as a thin push-off of $U$, and the second vector $e_2$ is perpendicular to $e_1$ in $\br^3$.  Or, analytically, taking $D$ to be the unit disk in the $r\theta$-plane in a coordinate chart $\br^4=\br^3\times\br$ (using cylindrical coordinates $(r,\theta,z)$ in $\br^3$ and coordinate $t$ in $\br$) the $0$-framing is $(\partial_z,\partial_r,\partial_t)$ and $1$-framing is $(\cos\theta\,\partial_z,\sin\theta\,\partial_r,\partial_t)$.
\end{comment}  
%%%%%%%%%%

To complete the proof of Theorem 1, it remains to show that the element $e_\bx$ in $\ext(\hone,\bztwo)$ defined by the extension ({\bf e}) maps to the Stiefel-Whitney class $\wtwox$ under the monomorphism
$$
\delta:\ext(H_1(\bx),\bztwo) \lto H^2(\bx;\bztwo) 
$$
in the universal coefficient sequence for $\bx$.  

To see this, first recall the definition of $\delta$.  View $\ext(H_1(\bx),\bztwo)$ as $\coker(i^*)$, where $i:\b\to\z$ is the inclusion of singular $1$-boundaries to $1$-cycles in $X$.  Then $\delta$ maps the equivalence class of a functional $\phi:\b\to\bztwo$ to the cohomology class of the cocylce $\phi\partial$, where  $\partial :\c\to \b$ is the boundary map on chains in $\bx$, that is, $\delta[\phi] = [\phi\partial]$. 

Alternatively, $\delta$ can be viewed as the dual of the Bockstein homomorphism
$$
\beta : H_2(\bx;\bztwo) \lto \tor_2(\hone)
$$
for the coefficient sequence $0\to\bz\overset{\times2}\to\bz\to\bztwo\to0$.  To explain this, we appeal to:

\begin{remark}\label{rem:ext}
There is a classical method for describing extensions by 2-torsion groups $T$ due to Eilenberg and MacLane \cite[Theorem 26.5]{EilenbergMacLane}. Let $A$ and $B$ be abelian groups fitting into an extension
$$
0\lto T \overset{u}{\lto} A\overset{h}{\lto} B\lto 0.
$$
Given any $b\in\tor_2(B)$, choose $a\in A$ with $h(a)=b$. Then there exists $t\in T$ such that $u(t) = 2a$, and it is easy to check that $t$ is uniquely determined by $b$ because $T$ consists of 2-torsion only. In fact, this leads to a homomorphism of groups $\tor_2(B) \to T$, inducing an isomorphism $\ext(B,T) \cong \hom(\tor_2(B),T)$.  
\end{remark}

For the case at hand we have $\ext(\hone,\bztwo)\cong\tor_2(\hone)^*$, mapping $[\phi]$ to the functional $[z]\mapsto\phi(2z)$, where $z$ is any integral cycle representing a 2-torsion element in $\hone$.  From this  one can easily check that the diagram
$$
\begin{CD}
\ext(\hone,\bztwo) @>\delta>>  H^2(\bx;\bztwo) \\
@V{\cong}VV @VV{\cong}V \\
\tor_2(\hone)^* @>>{\beta^*}> H_2(\bx;\bztwo)^*
\end{CD}
$$  
commutes, which is the sense in which $\delta$ and $\beta$ are dual.

Now $e_\bx$ is represented by any functional $f:\b\to\bztwo$ for which the diagram
$$
\begin{CD}
0 @>>> \b @>i>> \z @>>> \hone @>>> 0 \\ 
@. @VfVV @VVgV @| @. \\ 
0 @>>> \bztwo @>>u_1> \bfone @>>h_1> \hone @>>> 0
\end{CD}
$$
commutes; see for example Spanier \cite[\S5.5.2]{Spanier}.  This requires an initial choice of a map $g$ making the right square commute, and then $f$ is forced.  A different choice of $g$ will only change $f$ by the restriction of a functional on $\z$, and so its equivalence class $e_\bx = [f]$ in $\ext(\hone,\bztwo)$ is well-defined.   

Under the identification $\ext(\hone,\bztwo)\cong\tor_2(\hone)^*$, the element $e_\bx$ corresponds to the ``characteristic functional" 
$$
w:\tor_2(\hone)\lto\bztwo
$$
whose kernel is subgroup $\sone$ of split classes, introduced in \secref{prelim}.
%(see page~\pageref{twisted}).  
Thus if $K$ is a knot of order 2 in $\hone$, then $w(K)=0$ or $1$ according to whether $K$ (endowed with either framing) has order $2$ or $4$ in $\bfone$, or equivalently, whether $2K$ (meaning two copies of $K$ with the same framing on each) is framed bordant to $U_0$ or to $U_1$.  Note that up to framed bordism, the framing on $2K$ does not depend on the choice of framing on $K$, since this choice is being doubled.    

%%%%%%%%%
\begin{comment}
Equivalently, $w(K) = 0$ or $1$ according to whether
$h_1^{-1}[K]\cong\bztwo\oplus\bztwo$ or $\bz_4$.
\end{comment}
%%%%%%%%%

Thus to prove $\delta(e_\bx) = w_2(\bx)$, we must show that the functional $w_2\in H_2(\bx;\bztwo)^*$ corresponding to $\wtwox$ is equal to the composition $w\circ \beta$.  But $\beta$ sends the class of a surface $F$ to the class of any curve $C$ that is {\em characteristic} in $F$, meaning Poincar\'e dual to $w_1(F)$,\foot{This is well known, and easily verified.  Note that $C$ does in fact represent an element in $\tor_2(\hone)$:  It is null-homologous when $F$ is orientable, and of order $2$ in $H_1(F)$ otherwise, since $F$ can then be viewed as an integral $2$-chain with boundary $2C$.}  and $w_2$ reports self-intersections, by the Wu formula, so it remains to show
\begin{equation}\tag{{\bf w}}\label{eqn:w}
w(C) = \ F\dt F \pmod2
\end{equation}
for any surface $F\subset\bx$ and characteristic curve $C$ in $F$ .  

\begin{figure}[h!]
\includegraphics[height=140pt]{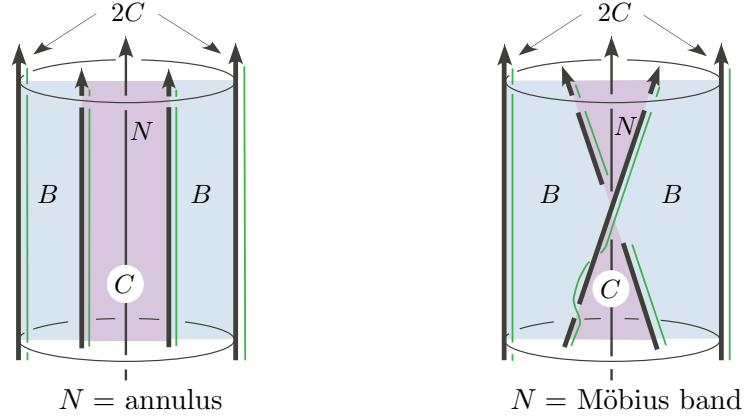}
\put(-235,137){\small$2C$}
\put(-234,34){\small$C$}
%\put(-233,82){\small$C$}
\put(-228,92){\small$N$}
%\put(-242,35){\small$N$}
\put(-205,68){\small$B$}
\put(-263,68){\small$B$}
\put(-51,137){\small$2C$}
\put(-50,32){\small$C$}
%\put(-49,90){\small$C$}
\put(-45,94){\small$N$}
%\put(-55,29){\small$N$}
\put(-27,68){\small$B$}
\put(-73,68){\small$B$}
\put(-255,-10){$N=$ annulus}
\put(-84,-10){$N=$ M\"obius band}
\caption{The bordism $B$ between $2C$ and $\partial N$}
\label{fig:bordism}
\end{figure}

If $F$ is orientable, both sides of equation ({\bf w}) vanish since $C$ is empty and $\bx$ is even. If $F$ is non-orientable, we can choose $C$ to be connected (it is then characterized up to homology by the orientability of its complement $F-C$) with closed tubular neighborhood $N$ in $F$.  Then $N$ is either an annulus or a M\"obius band, according to whether $C$ is orientation preserving or orientation reversing in $F$.   In either case, $\partial N$ is bordant in $X$ to $2C$; the bordism $B$ is trivial when $N$ is an annulus, and a pair of pants when $N$ is a M\"obius band, as indicated in \figref{bordism} with one suppressed dimension.  

%%% FIGURE 2: Bordism %%%
%%%%%%%%%%%%%%%%%%     

Starting with any framing on $C$ in $\bx$, the induced framing on $2C$ extends across $B$ (tilted in $\bx\times I$) to give a normal framing $e=(e_1,e_2,e_3)$ on $\partial N$.  We can assume that $e_1$ is normal to $N$ along $\partial N$, pointing in the suppressed direction in the figure, $e_2$ is shown there by the thin 
%(green) 
push-off, and $e_3$ is determined by $e_1$, $e_2$ and the orientation on $\partial N$.  Then $e_1$ extends to a normal vector field across all of $N$, and since $F-N$ is orientable, across the rest of $F$ off of a disk $D$; we use this vector field to push $F-D$ off itself in order to compute $F\dt F$.   Similarly, there is no obstruction to extending the framing $e$ on $\partial N$ across $F-(N\cup D)$ (tilted in $\bx\times I$) to give a framing on the unknot $U = \partial D$, which we also call $e$.  By construction $2C$ is framed bordant to $U_e$.  
Now the mod $2$ self-intersection $F\dt F$ is $0$ or $1$ according to whether the framing $e$ on $U = \partial D$ extends over $D$ or not, or equivalently, whether $U_e$ is framed bordant to $U_0$ or $U_1$.  Thus by definition $F\dt F= w(C)$, proving equation ({\bf w}) and completing the proof of Theorem 1.  
\qed

%%%%%%%%%%%%%     EXAMPLES     %%%%%%%%%%%%%%%%%%%
\begin{example}
If $\hone$ has no $2$-torsion, then it is immediate from the theorem that $\bfone$ is isomorphic to $\hone$ or to $\hone\oplus\bztwo$, depending upon whether $\bx$ is odd or even.  In particular for $\bx$ simply connected, $\bfone$ is trivial when $\bx$ is odd, and isomorphic to $\bz_2$ when it is even.

If $\hone$ has $2$-torsion, then 
$
\hone \ \cong \ \bz_{2^{k_1}}\oplus \cdots \oplus \bz_{2^{k_n}} \oplus A
$
for some sequence $1\le k_1\le \cdots \le k_n$ of integers, where $A$ has no $2$-torsion.  For each $i = 1,\dots,n$ let $C_i$ be a knot in $\bx$ representing $2^{k_i-1}$ times the generator of the $i^\text{th}$ summand above.  Then $C_1,\dots,C_n$ generate $\tor_2(\hone)$.   Now recall the functional $w:\tor_2(\hone)\longrightarrow\bztwo$ from  the proof of the theorem, given by $w(C) = F\dt F \pmod2$, where $F$ is any surface in $\bx$ in which $C$ is characteristic.   If $w\equiv0$, or equivalently $\bx$ is spin, then of course $\bfone \cong \hone\oplus\bztwo$.  Otherwise choose the smallest $i$ for which $w(C_i)=1$, and rewrite $\hone \cong \bz_{2^{k_i}} \oplus B$, where $B$ incorporates the remaining summands.  Then   
$$
\bfone \ \cong \ \bz_{2^{k_i+1}} \oplus B.
$$
This is a consequence of the identity $w\circ \beta=w_2$ established in the proof of the theorem.

From this observation, it is easy to construct $4$-manifolds $X$ with $\hone$ arbitrary, and with $\bfone$ equal to any prescribed $\bztwo$-extension of $\hone$.  For example, to get $\hone
 \cong \bz_{2^k}\oplus B$ and $\bfone \cong \bz_{2^{k+1}}\oplus B$, take $\bx$ to be the connected sum of the even manifold $E_{2^k,1}$ with suitable $E_{n,0}$'s (see \figref{even}).  The special case $\bff_1(E_{2,1}) = \bztwo$ and $\bff_1(E_{2,1}) = \bzfour$ was noted in Example 6.4 of Taylor\,\cite{Taylor}.   
\end{example}

%%%%%%%%%%%%%  REMARK   %%%%%%%%%%%%%%%%%
\begin{remark}\label{rem:Taylor}
We explain here how Theorem 1 can be deduced from Taylor's Theorem 6.2 in \cite{Taylor}.  As noted on page \pageref{prop:claim} above, Steenrod\,\cite[Theorem 28.1]{Steenrod} constructed an extension
$$
0 \lto T \lto [X,S^3]\lto H^3(X) \to 0
$$
for any 4-complex $X$, where $T=\coker(\,\overline{Sq}^2:H^2(\bx) \to H^4(\bx;\bztwo)\,)$; see Taylor \cite[6.1]{Taylor}.  By Remark~\ref{rem:ext}, this extension is classified by a homomorphism
$
\varepsilon_X : \tor_2(H^3(X)) \to T. 
$
The Bockstein homomorphism induces an isomorphism 
$$
\beta:H^2(X;\bztwo)/\im(H^2(X)) \lto \tor_2(H^3(X))
$$
and the Steenrod square $Sq^2: H^2(X;\bztwo) \to H^4(X;\bztwo) $ induces a homomorphism 
$$
[Sq^2]: H^2(X;\bztwo)/\im(H^2(X))\lto T.
$$
With a little bit of work, the statement of Theorem 6.2 in \cite{Taylor} can be interpreted as saying that $\varepsilon_X=  [Sq^2]\circ \beta^{-1}$. Since $Sq^2(x) = x\cup x$, this statement implies our Theorem 1.  

More directly, the dual of Taylor's Theorem 6.2 for the case $n=3$ and $k=1$ is exactly the assertion ({\bf w}), 
% on page \pageref{eqn:w}, 
which is the key step in our proof of Theorem 1. \end{remark}

\medskip
%%%%%%%%%%%%%%%%%%%%%%%%%%%%%%%%% %%%%%%%
%%    SEC 2:  ENUMERATION OF THE SECOND FRAMED BORDISM SET    %%
%%%%%%%%%%%%%%%%%%%%%%%%%%%%%%%%%%%%%%%% 
\section{Computation of the set $\bftwo \ \cong \ \pitwo$} \label{sec:two}
%%%%%%%%%%%%%%%%%%%%%%%%%%%%%%%%%%%%%%%%
%%%%%%%%%%%%%%%%%%%%%%%%%%%%%%%%%%%%%%%%

\begin{figure}[h!]
\hskip10pt \includegraphics[height=100pt]{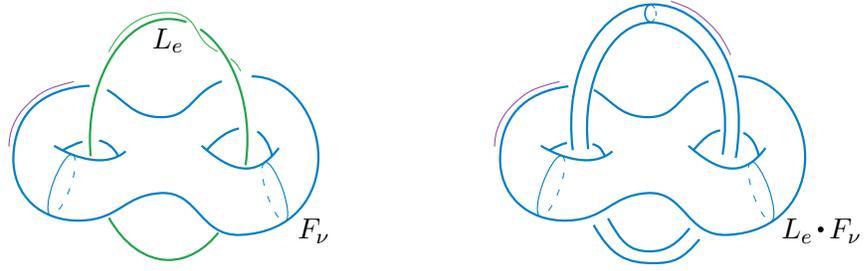}
\put(-248,82){$L_e$}
\put(-193,10){$F_\nu$}
\put(-10,10){$L_e\dt F_\nu$}
\caption{The framed link translate of a framed surface}
\label{fig:action}
\end{figure}

The forgetful map
$$
h_2 : \bftwo \longrightarrow \htwo
$$
sends the bordism class of a framed surface in $\bx$ to the homology class of the underlying oriented surface.  The image of $h_2$ is the set $\htwoo$ of all classes in $\htwo$ of self-intersection zero, that is, those represented by surfaces in $\bx$ with trivial normal bundle.  Thus to compute $\bftwo$, it suffices for each $\alpha\in\htwoo$ to classify up to framed bordism the framings on surfaces representing $\alpha$, that is, to enumerate the framed bordism classes in the fiber 
$$
\bfa \ = \ h_2^{-1}(\alpha)\ \subset \ \bftwo\,.
$$
To accomplish this, {\it fix a surface} $F$ representing $\alpha$ and {\it choose a framing} $\nu$ of the normal bundle of $F$.  Then $F_\nu$ represents one element in $\bfa$, and all others arise as ``framed link translates" of $F_\nu$, in the following sense:

Let $L$ be an oriented link in $\bx$ with normal framing $e = (e_1,e_2,e_3)$.  We can assume that $L$ is disjoint from $F$, after an isotopy if necessary.   Now form a new framed surface $L_e\dt F_\nu$, called the {\it $L_e$-translate of $F_\nu$}, by adjoining one new framed torus $T_{\tau(e)}$ for each component $K$ of $L$, where $T$ is the boundary of a small $3$-dimensional solid torus thickening $V$ of $K$ in $\bx$ in the direction of $e_1$ and $e_2$.  This is illustrated schematically in \figref{action}.

%%% FIGURE 3: Action %%%%%%%%%%
%%%%%%%%%%%%%%%%%%%%%%%

\begin{figure}[h!]
\hskip10pt \includegraphics[height=110pt]{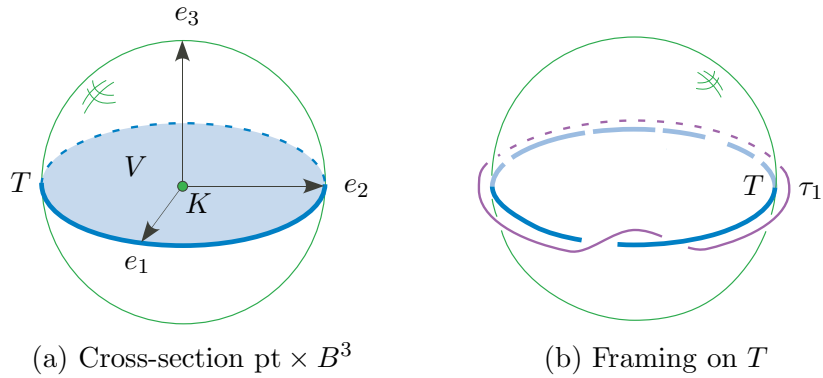}
\put(-260,22){$e_1$}
\put(-177,51){$e_2$}
\put(-241,116){$e_3$}
\put(-260,57){$V$}
\put(-303,51){$T$}
\put(-237,43){$K$}
\put(-26,50){$T$}
\put(-5,50){$\we_1$}
\put(-295,-15){(a) Cross-section $\pt\times B^3$\hskip 1in (b) Framing on $T$}
\caption{The framed torus $T=\partial V$ associated with the framed knot $K$}
\label{fig:framing}
\end{figure}

To describe the framing $\tau(e) = (\we_1,\we_2)$ on $T$ precisely, we set up coordinates as follows.  Identify a small tubular neighborhood of $K$ in $\bx$ with $S^1\times B^3$ via the framing $e$, with $K=S^1\times0$.  Then set 
$$
V \ = \ S^1\times D  \qquad\text{and}\qquad  T \ = \ \partial V \ = \ S^1\times \partial D
$$
where $D$ is the equatorial disk $B^3\,\cap\,\langle e_1,e_2\rangle$, as shown in \figref{framing}(a) at a crossection $\pt\times B^3$.  Now we require $\tau$ to spin once relative to the ``constant" framing $(\partialr,\,\partialz)$, using cylindrical coordinates $(r,\theta,z)$ in $B^3$, so that in particular it does {\it not} extend across the disk fibers of $V$ pushed into $\bx\times I$.  Explicitly $\we_1 = \partialr\cos\theta + \partialz\sin\theta$ and $\we_2=-\partialr\sin\theta + \partialz\cos\theta$ at any $\pt\times(1,\theta,0)\,\in\,T$, as illustrated in \figref{framing}(b) by showing (up to homotopy) the tip of $\we_1$ in each frame.

%%% FIGURE 4: Framing %%%
%%%%%%%%%%%%%%%%%%%%%%%

%%%%%%%%%%%%%     REMARK     %%%%%%%%%%%%%%%%%%%
\begin{remark}\label{rem:trivial}
The $0$-framed unknot $U_0$ acts trivially on any framed surface $F_\nu$, that is $U_0\dt F_\nu$ is framed bordant to $F_\nu$.  We leave this as an instructive exercise for the reader.  Of course in general, $L_e\dt F_\nu$ need not be framed bordant to $F_\nu$.  For example, if $F_\nu$ is a framed $2$-sphere in $S^4$ (note that all framings on embedded $2$-spheres in $4$-manifolds are isotopic since $\pi_2SO(2)=0$) then $U_1\dt F_\nu$ is not framed bordant to $F_\nu$.  Indeed $[U_1\dt F_\nu]$ generates $\bff_2(S^4) = \pi^2(S^4) = \pi_4(S^2) = \bztwo$, while $[F_\nu]=0$.
\end{remark}
%%%%%%%%%%%%%%%%%%%%%%%%%%%%%%%%%%%%%%%%%

There is an another more direct way to modify the framing $\nu$ using a link $L$ lying {\it on} $F$.  Simply add a full right handed twist in the normal fibers along any transverse arc to $L$ in $F$.  This results in a new framing on the same surface, which we call the {\it $L$-twist} of $\nu$ and, by abuse of notation, denote by $\nu+L$.  

In fact any framing $\hat\nu$ on $F$ can be obtained in this way.  Indeed $\hat\nu=\nu+L$ where $L$ is any embedded representative of the Poincar\'e dual of the cohomology class in $H^1(F)$ (with integer coefficients) that measures the difference $\hat\nu-\nu$ (meaning that this difference on any given $1$-cycle in $F$ is $k$ full right handed twists, where $k$ is the intersection number of the cycle with $L$).  Thus any two framings on $F$ are related by ``link twists".

It will be seen in the proof of the ``action lemma" below that $F_{\nu+L}$ is framed bordant to $L_e\dt F_\nu$, where $e$ is the framing on $L$ given by  $\nu$ together with the normal to $L$ in $F$, and so link twisting can be thought of as a special case of framed link translation.
  
%%%%%%%%%%%%%     ACTION LEMMA     %%%%%%%%%%%%%%%%%
\begin{actionlemma}
Framed link translation defines an action of the group $\bfone$ on the set $\bftwo$, given by $[L_e]\dt [F_\nu] =[L_e\dt F_\nu]$.  The orbit of $[F_\nu]$ under this action is the set $\bfa$ of all classes of framed surfaces that, forgetting the framing,  represent $\alpha=[F]\in\htwo$.
\end{actionlemma}  
%%%%%%%%%%%%%%%%%%%%%%%%%%%%%%%%%%%%%%%%%

\begin{proof}  Since addition in $\bfone$ is disjoint union, it is clear that the operation $\dt$ in the lemma defines an action {\sl provided} it is well-defined.  To see that it is in fact well-defined, assume that $L_e$ and $\hat L_{\hat e}$ are bordant via a framed surface $E^2\subset \bx\times I$, and that $F_\nu$ and $\hat F_{\hat\nu}$ are bordant via a framed $3$-manifold $Y^3\subset\bx\times I$.  We must show that $L_e\dt F_\nu$ and $\hat L_{\hat e}\dt \hat F_{\hat\nu}$ are framed bordant.  

To do so, first adjust $E$ and $Y$ to intersect transversally in finitely many points $x_1,\dots,x_n$.   Then ``pipe" these intersections to the upper boundary $\hat F$ of $Y$ in the usual way:  Choose disjoint arcs $\gamma_i$ in $Y$ joining the points $x_i$ to points $y_i$ in $\hat F$, and thicken these into disjoint solid cylinders $\gamma_i\times D$ in $Y$ that meet $E$ and $X\times 1$ in disks $E_i =  x_i\times B^2$ and  $D_i = y_i\times B^2$; the $D_i$ are meridional disks for $\hat F$ in $\bx\times 1$.  Now construct a new surface $\hat E$, removing the disks $E_i$ from $E$ and replacing them with the cylinders $\gamma_i\times \partial B^2$, and then rounding corners.  The framing on $E$ clearly extends across $\hat E$, restricting to the $0$-framing on each unknot $\partial D_i$.  Thus $\hat E$ provides a framed bordism between $L_e\dt F_\nu$ and $(\hat L_{\hat e}+L_0)\dt \hat F_{\hat\nu}$, where $L_0$ is a $0$-framed unlink in $\bx$.  But as noted in \remref{trivial} above, $L_0$ acts trivially, and so $L_e\dt F_\nu$ and $\hat L_{\hat e}\dt \hat F_{\hat\nu}$ are framed bordant.

To complete the proof of the action lemma, we must show that any framed surface $\hat F_{\hat\nu}$ belonging to $\bfa$ is framed bordant to $F_\nu$ acted on by some framed link $L_e$.   We can assume that $F$ and $\hat F$ are connected, and by hypothesis, they are bordant.  Hence we can find a bordism $Y$ from $F$ to $\hat F$ with only $1$ and $2$-handles.  The framing $\nu$ extends (not uniquely) over the $1$-handles to give a normal framing on the middle level of $Y$.  Similarly, $\hat\nu$ on $\hat F$ extends to a normal framing on the middle level.  Abusing notation, call the middle level $F$ and the two normal framings $\nu$ and $\hat\nu$.

Now let $L$ be a link in $F$ representing the Poincar\'e dual of the cohomology class in $H^1(F)$ that measures the difference $\hat\nu-\nu$, and so $\hat\nu=\nu+L$ in the notation above.   Then $L$ acquires a framing $e$ by appending its normal vector field in $F$ onto $\nu$, and we propose to show that $F_{\nu+L}$ is framed bordant to $L_e\dt F_\nu$.

To that end, we first ``blister" $L$. That is, we create a $3$-dimensional bordism $Y^3$ obtained from $F \times I$ (using the first vector of $\nu$) by removing an open disk bundle neighborhood of $L \times 1/2$ so that $\partial Y = (F \times 0) \cup (F \times 1) \cup T$ where $T$ is the torus boundary of the disk bundle.  Then ``tilt" $Y$ into $X \times I$ so that $Y\cap (X \times 0) = (F\times 0)\cup T$ and $Y \cap (X \times 1)= F \times 1$. 

%%% FIGURE 5: Blister %%%%%%%%%%%
\begin{figure}[h!]
\hskip10pt \includegraphics[height=110pt]{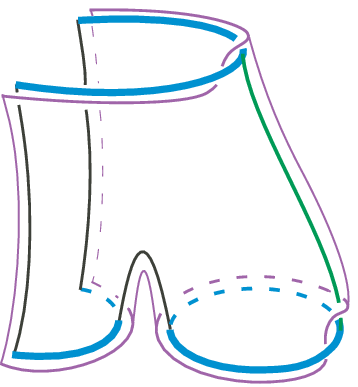}
\put(-90,-10){$(F\times0)_\nu\cup T_{\tau(e)}$}
\put(-85,117){$(F\times1)_{\hat\nu}$}
\put(-33,70){\small$\alpha$}
\put(-40,45){$Y^3$}
\caption{Blistering the link $L$}
\label{fig:blister}
\end{figure}
%%%%%%%%%%%%%%%%%%%%%%%

The framing $e$ on $L$ induces a framing $\tau(e)$ on $T$ by the link translation construction above.  Together these give a framing $\nu\cup\tau(e)$ on $(F\times 0)\cup T$ that extends over $Y$ by adding an extra full twist when going past the arc $\alpha = Y\cap(L \times [1/2,1])$, the sign of the twist depending on the orientations of $L$ and $F$.  By construction, this framing restricts to $\hat\nu = \nu+L$ on $F \times 1$.  This is shown in \figref{blister} with two of the five dimensions suppressed, one tangent to $F$ and the other normal to $F$ in $\bx$.
Thus we have constructed a framed bordism $Y$ between $F_\nu\cup T_{\tau(e)}= (F \times 0)_\nu \cup T_{\tau(e)}$ and $F_{\hat\nu} = (F \times 1)_ {\hat\nu}$, and the former is the action of $L_e$ on $F_\nu$.  This completes the proof.
\end{proof}

By the action lemma, the classes in $\bfa=h_2^{-1}(\alpha)$ correspond to the cosets in $\bfone$ of the stabilizer $\stab$ of any chosen $\alphat= [F_\nu]$ in $\bfa$ under the $\bfone$-action.  We will show that this stabilizer is exactly  {\it twice} the image of the {\it framed intersection map}
$$
i_\alphat : \bfthree \lto \bfone
$$
that sends (the class of) an oriented $3$-manifold $Y$ in $\bx$ with normal $n$ to its transverse intersection with $F$, framed by $(\nu,n)$.  Note that $n$ is determined (up to homotopy) by the orientation on $Y$, using the convention that $(\tau,n)$ should frame $T\bx$ over $Y$ for any oriented framing $\tau$ of $TY$, and so $i_\alphat$ can be viewed as a map $\hthree\to\bfone$.  A straightforward transversality argument shows that $i_\alphat$ is a well defined homomorphism.

Since there are only two framings on any 1-dimensional link in $\bx$, up to framed bordism, the map $2i_\alphat$ depends only on $\alpha$ (independent of the choice of $\alphat$) and so we denote it by 
$$
\tau_\alpha\ = \ 2\,i_\alphat:\bfthree\lto\bfone\,.
$$
Thus we claim that $\stab = \im(\tau_\alpha)$.  This will complete the proof of Theorem 2, restated below.  The first part of the theorem is just a summary of the discussion above, with the action lemma providing the proof that the $\bfone$-action on $\bfa$ is transitive. 

%%%%%%%%%%%%%     THEOREM 2     %%%%%%%%%%%%%%%%
\begin{theorem2}
The forgetful map $h_2:\bftwo \to H_2(\bx)$ has image equal to the set $\htwoo$ of all classes of self-intersection zero.  The fiber $\bfa = h_2^{-1}(\alpha)$ over any $\alpha\in\htwoo$ has a transitive $\bfone$-action by framed link translation {\rm (defined above)} with stabilizer equal to the image of the map $\tau_\alpha:\bfthree\to\bfone$ that records {\rm\bf twice} the intersection with any framed surface representing an element in $\bfa$.  Thus $\bfa$ is in (non-canonical) one-to-one correspondence with the cokernel \,$\fr$ of $\tau_\alpha$.
\end{theorem2}%%%%%%%%%%%%%%%%%%%%%%%%%%%%%%%%%%%%%%

\begin{proof}  Given a class $\alpha\in\htwoo$, fix an oriented surface $F$ in $X$ representing $\alpha$ and a framing $\nu = (n_1,n_2)$ of its normal bundle.  Set $\alphat = [F_\nu]\in \bfa$ and let $\stab$ denote the stabilizer of $\alphat$ under the $\bfone$-action on $\bfa$.   It remains to prove  
$$
\stab \ = \ \im(2\,i_\alphat).
$$

We first verify the inclusion 
$$
\im(2\,i_\alphat) \ \subset \ \stab.
$$
Let $Y \subset \bx$ be an oriented $3$-manifold transverse to $F$, and $L = Y \cap F$ with framing $e = (n_1,n_2,n)$, where $n$ is the oriented normal to $Y$ in $\bx$. Thus $L_e$ represents a typical element of $\im(i_\alphat)$.  Now identify $\bx$ with $\bx\times1/2$.  Then $Y\subset \bx\times I$ has a natural framing $(n,t)$, where $t$ is tangent to the $I$-factor, as does $F\times I$, given by $\nu$ at all levels.  The union $Y \cup (F \times I)$ is a framed, {\it immersed} submanifold of $X\times I$ which has, using the framing $\nu$ on $F \times\partial I$, self-intersection equal to two copies of $L$.  We propose to replace this with an {\it embedded} submanifold, keeping track of how the framing changes.

This is analogous to the well known case of a $2$-sphere $S$ immersed in the $4$-sphere with one positively oriented double point $p$.  The self intersection $S\dt S$ must be zero since $S$ can be pushed far off of itself.  But $S\dt S$ is also given by pushing $S$ off itself using its normal bundle and counting points of intersection;  the double point $p$ provides two plus points, so there must be two other minus points.  

An alternative description is an immersed $2$-disk with one double point where the same push off as above leads to a framing on the bounding circle which has two full negative twists compared to the trivial framing.  On the other hand, the double point can be removed by replacing the two transverse disks by an annulus, obtaining an embedded torus or punctured torus.  The torus has a trivial normal bundle, so the punctured torus gets the trivial framing on the normal bundle of its bounding circle.  Thus, removing the double point via the annulus changes the framing by two positive full twists.

Similarly, removing the double points using annuli all along the points of $L = Y \cap (F \times I)$ changes the framing on (say) $F\times 1$ by $2L$, meaning that the new framing on $F\times 1$ is $\nu+2L$, the ``$2L$-twist of $\nu$" in the terminology introduced above.  Thus the two framed surfaces $F _ \nu$ and $F_{ \nu + 2L}$ are framed bordant  by $Y \cup (F \times I)$ made into an embedding by removing the double points, as above.  The proof of the action lemma then shows that $F_{\nu + 2L}$ is framed bordant to $2L_e\dt F_\nu$, and so $2L_e$ stabilizes $F_\nu$.

For the reverse inclusion 
$$
\stab \ \subset \ \im(2\,i_\alphat)
$$
start with a framed link $L_e$ representing an element in $\stab$.  This means there is a framed bordism $W_\mu\subset \bx\times I$ from $F_\nu$ to $\hat F_{\hat\nu} =L_e\dt F_\nu$.  As in the proof of the action lemma, we can assume that $F=\hat F$, that $L$ lies in $F$ where it is dual to the class in $H^1(F;\bz)$ that measures the difference between $\hat\nu$ and $\nu$, and that the framing $e$ on $L$ is induced from $\nu$.   We must find a closed $3$-manifold $Y\subset\bx$ such that 
$$
2\,i_\alphat(Y) \ = \ 2[Y\cap F] \ = \ [L_e] \ \in \ \bfone.
$$ 

To construct $Y$, start with the $3$-cycle $V = W - (F \times I)$.  Since $W$ and $F\times I$ have the same boundary, this cycle is homologous to an embedded $3$-manifold in $X = X\times 1/2$, which can be assumed transverse to $F$.  We let $Y$ be any such $3$-manifold, and set $J=Y\cap F$.  

It remains to show that $2J$ is homologous to $L$ in $\bx$, or equivalently in $\bx\times I$.   Working in $H_1(\bx\times I)$, it is clear that $V^2 = V\dt V=0$, since $V$ is homologous to a submanifold of $\bx\times1/2$.  But this self-intersection can also be computed in terms of intersections of the relative cycles $W$ and $F\times I$, as follows:
$$
\begin{aligned}
V^2 \ &= \ (W-(F\times I))^2 \\
&= \ W^2 - 2\,W\dt(F\times I) +  (F\times I)^2
\end{aligned}
$$
Here it is understood that when computing intersection numbers in $\bx\times I$,  we always use $\nu$ and $\hat\nu$ to start the push off along the boundary.  With this understanding we have $W^2 = 0$, and so 
$$
2\,W\dt(F\times I) \ = \ (F\times I)^2.
$$
Now it is clear that the left hand side is represented by $2J$, while the right hand side is represented by $L$.  Hence $2J$ and $L$ are homologous, and the theorem is proved.  
\end{proof}

%%%%%%%%%%%%%     REMARK     %%%%%%%%%%%%%%%%%
\begin{remark}\label{rem:coker}
The group $\fr = \coker(\tau_\alpha)$ featured in Theorem 2 can be computed in terms of the group $\ufr = \coker(2\,i_\alpha)$, where $i_\alpha$ is the (unframed) intersection map
$$
i_\alpha:\hthree\longrightarrow\hone
\qquad 
\beta\longmapsto \alpha\dt\beta.
$$
{\em For example \,$\fr \,\cong\, \ufr$\, when $\bx$ is odd}, since the forgetful map $h_1$ is then an isomorphism.  

If $\bx$ is even, then the sequence $({\bf e})$ in Theorem 1 induces a commutative diagram
$$
\begin{CD}
@.  @. \bfthree @>\cong>> \hthree @. \\ 
@. @. @V\tau_\alpha V{=\ 2i_\alphat}V @VV2i_\alpha V@. \\ 
0 @>>> \bztwo @>u_1>> \bfone @>h_1>> \hone @>>> 0 \\
@. @V{\cong}VV @VVV @VVpV@. \\ 
@. \bztwo @>u_1^\alpha>> \fr @>h_1^\alpha>> \ufr @>>> 0 \\
@. @. @VVV @VVV@. \\ 
@.  @. 0 @. 0 @. \\
\end{CD}
$$
with exact rows and columns.  Note that $u_1^\alpha$ is zero precisely when $u_1(1) = [U_1]\in\im(\tau_\alpha)$ (recall that $U_1$ is the $1$-framed unknot), or equivalently, when $\alpha$ is a twisted class (see \secref{prelim}). %page \pageref{twisted}).  
{\em Hence $\fr \,\cong\, \ufr$ when $\bx$ is of type III$_1$ and $\alpha$ is twisted}.

{\em In all other cases} \,($\bx$ of type II or type III$_2$, when there are no twisted classes, or of type III$_1$ with $\alpha$ untwisted)\, {\em the group $\fr$ is a $\bztwo$-extension of $\ufr$}, via the exact sequence
$$
0 \,\lto\, \bztwo \overset{u_1^\alpha}{\,\lto\,} \fr \overset{h_1^\alpha}{\,\lto\,} \ufr \,\lto\, 0\,.
$$
By \remref{ext} this extension is classified by the homomorphism $\tor_2(\ufr) \to \bztwo$ with kernel $p(\sone)$, where $\sone\subset\hone$ is the subgroup of split classes.  It follows that {\em if $\bx$ is of type II, then the extension splits {\rm (since $\tone=\O$)}, so \,$\fr \,\cong\, \ufr\oplus\bztwo$}.  This can also be seen directly from the diagram above, since $h_1$ splits by Theorem 1, and so $h_1^\alpha$ splits as well.  In contrast, {\em if $\bx$ is of type III and $\alpha$ is untwisted, then the extension does not split}, since $\im(2i_\alpha)$ will then contain some twisted $1$-dimensional classes.    

\end{remark}
%%%%%%%%%%%%%%%%%%%%%%%%%%%%%%%%%%%%%

%%%%%%%%%%%%  EXAMPLES   %%%%%%%%%%%%%

%%% EXAMPLE 1: Simply connected 4-manifolds %%%

\begin{example}
Let $\bx$ be simply connected.  Then $\hone$ and $\hthree = 0$.  

If (the intersection form on) $\bx$ is definite, then $\htwoo=\{0\}$, that is, {\it only} the trivial homology class has self-intersection zero.   Therefore 
$$
\bftwo \ \cong\  \bfone \ =\  0 \text{ \ or \ } \bztwo
$$
according to whether $\bx$ is odd or even ($=$ spin, since $\bx$ is simply connected).  Here the symbol $\cong$ indicates {\it torsor equivalence}, so in the even case there is no {\it a priori} way to distinguish between the two homotopy classes of maps $\bx\to S^2$.   In this case, however, it is natural to associate the constant map with $0\in\bztwo$ (as is the case more generally for $\bff_2^0 = h_2^{-1}(0) \subset \bftwo$).

If $\bx$ is indefinite and odd, then it is homotopy equivalent to $p\,\bc P^2\,\#\,q\,\overline{\bc P}^2$ for some $p,q>0$.  In this case $\htwoo$ is infinite, but each of its elements has a {\it unique} associated homotopy class of maps $\bx\to S^2$ since $\bfone = 0$.  Thus 
$$
\bftwo \ = \ \htwoo.
$$
The problem of enumerating the classes in $\htwoo$ is a classical and generally intractable number theoretic problem; these classes correspond to the different ways that integers can be expressed simultaneously as the sum of $p$ squares and as the sum of $q$ squares. 

There is a similar story in the indefinite and even case, except that now each class in $\htwoo$ has {\it exactly two} corresponding homotopy classes of maps to the $2$-sphere, since $\bfone = \bztwo$.  Hence
$$
\bftwo \ \cong \ \htwoo \times \bztwo
$$
as sheaves of torsors.  This means that there is no natural way to distinguish between the two homotopy classes of maps $\bx\to S^2$ associated with each nonzero class in $\htwoo$.
\end{example}

%%% EXAMPLE 2: Products Y^3 \times S^1 %%%

\begin{example}  \ Let $\bx = Y\times S^1$, where $Y$ is a closed, oriented $3$-manifold.   By the K\"unneth  theorem, there is a natural isomorphism
$$
\hk \ \cong \ H_k(Y)\oplus H_{k-1}(Y)
$$
for all $k$, with integer coefficients understood as always.  Elements in the first summand correspond to $k\text{-manifolds}$ in $Y\times\pt$, while those in the second correspond to products of $(k-1)$-manifolds in $Y$ with the circle.  

In particular, we will identify elements $\alpha\in\htwo$ with pairs $(\beta,\gamma)\in H_2(Y) \oplus H_1(Y)$.  It follows that $\bx$ is spin, since $\htwo$ is generated by elements of the form $(\beta,0)$ and $(0,\gamma)$ which have self-intersection zero in $\bx$, even allowing $\bztwo$ coefficients.  

We wish to compute the fiber $\bfa = h_2^{-1}(\alpha)$ over an arbitrary element $\alpha = (\beta,\gamma)$ in $\htwoo = \{(\beta,\gamma) \st \beta\dt\gamma = 0 \text{ \,in $Y$}\}$.  By Theorem 2,  this fiber is an $\fr$-torsor, where 
$$
\fr \ = \  \coker(2i_\alpha) \oplus \bztwo
$$
by \remref{coker}, since $\bx$ is spin.  Thus it suffices to compute $\coker(2i_\alpha)$.

Recall that $i_\alpha:\hthree\to\hone$ is the map that computes intersection with $\alpha=(\beta,\gamma)$.  Identifying $\hthree \cong H_3(Y)\oplus H_2(Y)$ and $\hone \cong H_1(Y) \oplus H_0(Y)$ as above, we have
$$
i_\alpha: \bz\oplus H_2(Y) \lto H_1(Y)\oplus\bz
$$
sends $(1,0)$ to $(\gamma,0)$ and $(0,\delta)$ to $(\delta\dt\beta,\delta\dt\gamma)$.  Hence $\coker(2i_\alpha)$ is isomorphic to the quotient of $H_1(Y)\oplus\bz$ by the subgroup generated by $(2\gamma,0)$ and twice the image of the map 
$$
i_\beta\oplus i_\gamma:H_2(Y)\lto H_1(Y)\oplus\bz
$$
whose coordinates record intersections in $Y$ with $\beta$ and $\gamma$.  In particular
$$
\bff_2^{(\beta,0)} \ \cong \ \coker(2i_\beta)\oplus\bz\oplus\bztwo   \and    \bff_2^{(0,\gamma)} \ \cong \ H_1(Y)/\langle2\gamma\rangle \oplus \coker(2i_\gamma) \oplus \bztwo.
$$
Note that $\coker(2i_\gamma) = \bz_{2d(\gamma)}$, where $d(\gamma)$ is the {\it divisibility} of $\gamma$ (the largest integer $d$ for which $d\cdot \kappa = \gamma$ has a solution $\kappa\in H_1(Y)$, or zero if $\gamma$ is of finite order).

We illustrate these computations with three cases: (a) $Y$ is a lens space, (b) $Y = S^2\times S^1$, and (c) $Y$ is the $3$-torus:

{\bf (a)}  Let $\bx = L(p,q)\times S^1$.  Then $\htwoo = \htwo \cong 0\oplus H_1(L(p,q)) \cong \bz_p$, and 
$$
\bff_2^k \ \cong \ \bz_{\gcd(2k,p)}\oplus\bz\oplus\bztwo.
$$
for any $k\in\bz_p$.  This is an instance of the case $(0,\gamma)$ above in which $i_\gamma = 0$.  

{\bf (b)}  Let $\bx = S^2\times T^2 = (S^2\times S^1)\times S^1$.  Then $\htwo = \bztwo$, generated by $\sigma = [S^2\times\pt]$ and $\tau = [\pt\times T^2]$.  The set $\htwoo$ consists of all integral multiples of $\sigma$ or $\tau$, and
$$
\bff_2^{d\sigma} \ \cong \  \bz \,\oplus\, \bz \,\oplus\, \bztwo 
\qquad\text{and}\qquad 
\bff_2^{d\tau} \ \cong \  \bz_{2d} \,\oplus\, \bz_{2d} \,\oplus\, \bztwo 
$$
for any integer $d$ (cf.\ Example 6.11 in \cite{Taylor}).  

{\bf (c)}  Let $\bx = T^4 = T^3\times S^1$.  Then $H_2(T^4) \cong H_2(T^3)\oplus H_1(T^3) \cong \bz^3\oplus\bz^3$, where the $\bz^3$ factors have dual bases $\beta_{23}, \beta_{31}, \beta_{12}$ (represented by the coordinate tori in $T^3$) and $\gamma_1, \gamma_2, \gamma_3$ (represented by the coordinate circles in $T^3$).  The matrix $B$ for the intersection pairing $H_2(T^3)\otimes H_2(T^3)\to H_1(T^3)$ is given by
$$
B \ = \ \begin{pmatrix} 
0&\gamma_3&-\gamma_2 \\ 
-\gamma_3&0&\gamma_1 \\ 
\gamma_2&-\gamma_1&0
\end{pmatrix}
$$
with respect to the ordered basis $\beta_{23}, \beta_{31}, \beta_{12}$ of $H_2(T^3)$.   The set $H_2^\circ(T^4)$ consists of all pairs $(\beta,\gamma)\in\bz^3\oplus\bz^3$ with $\beta\dt\gamma = 0$, where $\beta$ and  $\gamma$ are the coordinate vectors for elements in $H_2(T^3)$ and $H_1(T^3)$ with respect to the bases above, and $\dt$ is the usual dot product.  It follows that $\bff_{(\beta,\gamma)} \cong \coker(2A)\oplus\bztwo$, where
$$
A \ = \ \begin{pmatrix} 
0 & b_3 & -b_2 & -c_1 \\ 
-b_3 & 0 & b_1 & -c_2 \\ 
b_2 & -b_1 & 0 & -c_3 \\ 
c_1 & c_2 & c_3 & 0
\end{pmatrix},
$$
for any $(\beta,\gamma) = ((b_1,b_2,b_3),(c_1,c_2,c_3))$ in $H_2^\circ(T^4)$.  

To compute $\coker(2A)$, we reduce $A$ to its Smith Normal Form.   Note that $A$ is skew-symmetric and of determinant zero, since $\beta\dt\gamma = 0$.  Thus either $A=0$ or $\rk(A) = 2$.   In the first case set $d=0$, and in the second set $d$ equal to the greatest common divisor of all the entries in $A$.  Then by a change of basis of the form $P^TAP$, where $P$ is a unimodular integral matrix, we can assume $b_3=d$.  A further reduction shows that the Smith Normal Form of $A$ is the diagonal matrix with entries $d,d,0,0$.  Thus $\coker(A) = \bz_{d}^2\oplus \bz^2$, and so 
$$
\bff_2^\alpha \ \cong \  \bz_{2d} \oplus \bz_{2d} \oplus \bz \oplus \bz \oplus \bztwo
$$
for all $\alpha$ in $H_2^\circ(T^4)$, where $d$ is the divisibility of $\alpha$ in $H_2(T^4)$ .
\end{example}

%%% EXAMPLE 1: Type III_1 4-manifold %%%

\begin{example}\label{ex:Taylor}

To conclude, we describe a 4-manifold of type III$_1$, that is,  an even non-spin 4-manifold $\bx$ that has at least one twisted $2$-dimensional homology class $\alpha$ (see \secref{prelim}).  By \remref{coker}, the associated fiber $\bff_\alpha = h_2^{-1}(\alpha)$ is in one-to-one correspondence with the cokernel of $2i_\alpha:\hthree\to\hone$.  This shows that the last case discussed in Remark 6.8 in Taylor \cite{Taylor} is indeed realized by a 4-manifold.   

We construct $\bx$ by surgering the connected sum $E_{2,1} \,\#\, (T^2 \times S^2)$ along the curve $C \,\#\, J$, where $C$ generates $H_1(E_{2,1}) = \bztwo$ and $J$ is an essential curve in $T^2\times\pt$.  A handlebody for $X$, minus the 3 and 4-handles, is shown in \figref{twisted}, where $C$ can be taken to be the core of the obvious M\"obius band $M$ bounded by the attaching circle of the $1$-framed 2-handle (on the left).

%%% FIGURE 6: Type III_1 %%%
\begin{figure}[h!]
\includegraphics[height=100pt]{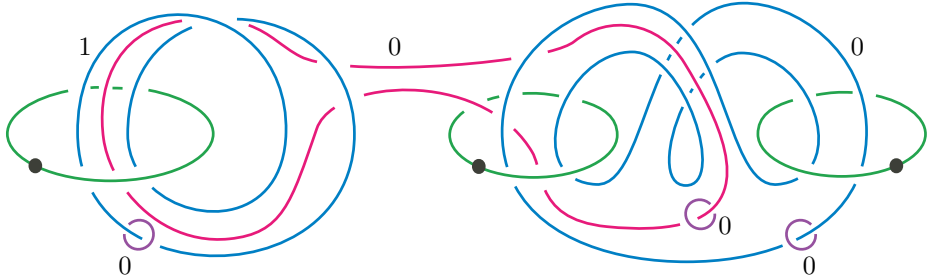}
\put(-322,80){\small 1}
\put(-307,-3){\small 0}
\put(-80,12){\small 0}
\put(-205,80){\small 0}
\put(-48,-3){\small 0}
\put(-30,80){\small 0}
\caption{A 4-manifold of type III$_1$}
\label{fig:twisted}
\end{figure}
%%%%%%%%%%%%%%%%%% 
\end{example}

\break

A straightforward computation gives
\begin{enumerate}
\item $\hone = \bztwo\oplus\bz$, with generators $\kappa_1 = [C]$ and $\kappa_2 = [K]$ where $K$ is a dual curve to $J$ in $T^2\times\pt$.
\item $\htwo = \bztwo\oplus\bz\oplus\bz$, with generators $\alpha_1$, represented by an $S^2$-fiber of the bundle $E_{2,1}$, $\alpha_2 = [T^2\times\pt]$ and $\alpha_3 = [\pt\times S^2]$ (in $T^2\times S^2$).
\item $\hthree = \bz$, with generator $\tau = [J\times S^2]$.
\end{enumerate} 

The intersection  pairing on $\htwo$ vanishes except for $\alpha_2\dt\alpha_3=1$, and so 
$$
\htwoo \ = \  \{(a,b,c)\in\bztwo\oplus\bz\oplus\bz \st bc = 0\}.
$$

The intersection pairing between $\htwo$ and $\hthree$ vanishes except that $\alpha_2\dt\tau = \kappa_1$, since $(T^2\times\pt)\cap(J\times S^2) = J$ which is isotopic to $C$ via the surgery.  Since $2\kappa_1 = 0$, the map $2i_\alpha$ is zero for all $\alpha\in\htwo$, and so $\coker(2i_\alpha) = \hone = \bztwo\oplus\bz$ for all $\alpha$.

Now observe that $\kappa_1=[C]$ is twisted.  Indeed $C$ is characteristic in an $\br P^2$-section of the bundle $E_{2,1}$ (the union of $M$ and the core of the $1$-framed 2-handle in the picture) and this section has odd self-intersection.  Hence $w(C)=1$ by equation ({\bf w}) in the proof of Theorem 1, where $w$ is the characteristic functional, and so $\kappa_1$ is twisted.  

It follows that $X$ has (infinitely many) $2$-dimensional twisted classes, namely all classes $\alpha = (a,b,c)$ in $\htwo$ with $b$ odd and $c=0$.   For each of these, $\bff_2^\alpha \cong \coker(2i_\alpha) = \bztwo\oplus\bz$.  For all other classes $\alpha$ in $\htwoo$ we have $\bff_2^\alpha \cong \bz_4\oplus\bz$, the nonsplit $\bztwo$-extension of $\hone$.

%%%%%%%%%%%%%%%%%%%%%%%%%%%%%%%%% %%%%%%%
%%%%%%%%%%%    SEC 3:  PROOF OF Theorem 3    %%%%%%%%%%%%
%%%%%%%%%%%%%%%%%%%%%%%%%%%%%%%%%%%%%%%% 
\section{Proof of Theorem 3} \label{sec:three}
%%%%%%%%%%%%%%%%%%%%%%%%%%%%%%%%%%%%%%%%
%%%%%%%%%%%%%%%%%%%%%%%%%%%%%%%%%%%%%%%%

For a topological group $G$, Milnor \cite{Milnor56} constructed a space $BG$, together with a principal $G$-bundle $EG\to BG$ such that the homotopy groups of $EG$ vanish. This implies by obstruction theory that isomorphism classes of $G$-principal bundles over a CW-complex $X$ are in 1-1 correspondence with $[X,BG]$. 

Given a closed subgroup $T<G$, we can form the associated bundle 
$$
EG \times_G G/T \lto EG/G
$$
with fiber $G/T$. The total space is homeomorphic to $EG/T$ and since $T$ is closed, the projection $EG \to EG/T$ is a principal $T$-bundle. By the weak contractibility of $EG$, the space $BT:= EG/T$ is a classifying space for principal $T$-bundles (over CW-complexes). 

Summarizing, we have a  fiber bundle 
$$
G/T \lto BT \overset{\pi}{\lto} BG
$$
as claimed in the introduction.  Since fiber bundles are Serre fibrations, the usual homotopy lifting theory gives an ``exact sequence'' for any CW-complex $X$:
$$
\begin{CD}
[X, \call BT]_{\bar u} @>\call\pi\,\scriptscriptstyle\circ\,_->> [X,\call BG]_b @>>> [X, G/T] @>h>> [X, BT] @>>> [X,BG]
\end{CD}
$$
Here we fix a map $ u: X\to G/T$ and consider its image $\bar u:X\to BT$, whose further image under $\pi$ is a constant $b\in BG$.  Exactness means that the group $[X,\call BG]_b$ acts transitively on the fiber $h^{-1}(\bar u)$ with stabilizer of $u$ the image of the map induced by $\call\pi$. Here $\call BT$ is the free loop space of $BT$  and the subscript $\bar u$ indicates that the maps start and end at $\bar u$, and similarly for $\call BG$. 

Encouraged by the referee, we add a short proof of this ``exact sequence'' for any Serre fibration $\pi:E{\to} B$ (in place of $BT\to BG$).  

\begin{proof}
Given a Serre fibration $E\overset{\pi}{\to} B$ with fibre $F$ over $b\in B$, there is  an ``exact sequence''
\[
\pi_1(E,e) \overset{\pi_*}{\ra} \pi_1(B,b) \ra \pi_0(F) \overset{h}{\ra} \pi_0(E) \ra \pi_0(B)
\]
in the following sense: For a fixed $e\in E$ with $\pi(e)=b$, the group $\pi_1(B,b)$ acts transitively on $h^{-1}(e)$, with stabilizer of $e$ equal to the image of $\pi_1(E,e)$ under $\pi_*$. This follows from the lifting properties of Serre fibrations, applied to intervals and squares. 

We next recall a useful tool, namely Steenrod's {\em convenient category $\CG$} of compactly generated Hausdorff spaces \cite{Steenrod2}. A topological space is {\em compactly generated} if a subset is closed provided its intersection with any compact set is closed. 

It is even more convenient to use Rainer Vogt's   (full) subcategory $\HG\subset \TOP$ which contains $\CG$. He shows in \cite[1.4]{Vogt} that the inclusion has a right adjoint $k:\TOP \to \HG$ where $k(X)$ is a finer topology on the underlying set of a topological space $X$. In fact, $k(X)$ is the finest topology such that any continuous map $K\to X$ from a compact Hausdorff space $K$ factors through the identity $k(X) \to X$. 

These categories are convenient for various reasons, one being that \cite[3.6]{Vogt} the canonical adjunction maps are homeomorphisms whenever $X,Y\in\HG$:
\[
B^{X \otimes Y} \ \cong \ (B^X)^Y
\]
Here $\otimes$ denotes the categorical product in $\HG$ and $B^X$ is the set of all continuous maps $X\to B$, equipped with the topology given by applying the functor $k$ to the compact open topology. It is important to recall from \cite[3.11]{Vogt} that for (locally) compact spaces $Y$, the categorical product $X \otimes Y$ is actually just the cartesian product $X \times Y$.

We note that any CW-complex $X$ is Hausdorff and compactly generated and therefore $X\in \CG\subset\HG$. Take our Serre fibration $\pi:E\to B$ and consider the induced map 
\[
\pi^X: E^X \ra B^X
\]
The crucial claim that this is again a Serre fibration: By the above adjunctions, lifting maps from a finite CW-complex $D \to B^X$ over $E^X$ is the same as lifting maps from $D \times X \to B$ over $E$. Since $D \times X$ is again a CW-complex, our crucial claim follows from the facts that 
\begin{itemize}
\item crossing with $X$ preserves acyclic cofibrations, and
\item it suffices to check the lifting conditions for Serre fibrations on {\em finite} CW-complexes.
\end{itemize} 

Now we apply our baby version of the ``exact sequence'' (for $\pi_1$ and $\pi_0$) to this new Serre fibration $\pi^X:E^X\to B^X$ and get our desired sequence:
\[
\begin{CD}
[X, \call E]_{\bar u} @>\call\pi\,\scriptscriptstyle\circ\,_->> [X,\call B]_b @>>> [X, F] @>h>> [X, E] @>>> [X,B]
\end{CD}
\]
The only thing to check are the following translations:
\[
\pi_0(B^X) \= B^X/\simeq \ \= [X, B]
\]
because a path $I\to B^X$ translates via our adjunction to a homotopy $X \times I \to B$ (since $I$ is compact Hausdorff). Moreover, if $b:X\to B$ denotes the constant map at $b\in B$ then 
\[
\pi_1(B^X,b) \= (B^X,b)^{(S^1,s)} / \simeq\ \= (B,b)^{(S^1 \times X, s \times X)} / \simeq\ \= [X, \call B]_b 
\]
where the right hand side denotes maps $X\to \call B$ that start and end at $b$. It would be natural to use the topology on $\call B$ which is given by applying the functor $k$ to the compact open topology on the space of maps $S^1\to B$. Note however, that the identity $k(Y)\to Y$ is always a weak equivalence and because $X$ is a CW-complex, the set $[X,\call B]$ does not see the difference between these two topologies.
\end{proof}

Going through the above arguments, we see that the action of $[X,\call BG]_b$ on $h^{-1}(\bar u)$ is given by first lifting a homotopy $H: X \times [0,1] \to BG$ that starts and ends at $b$,\foot{By adjunction, elements in $[X,\call BG]_b$ are represented by such homotopies, which are stacked to obtain the group structure, and similarly for $[X,\call BT]_{\bar u}$.} with prescribed lift $\bar u$ at time 0.  Then one evaluates the lift at time 1 to get a second map $X\to BT$ lying over the constant $b$ which can be identified with a map $X\to G/T$. 

Let $\Omega_b$ denote the space of loops based at $b$ and $\calp_b$ the space of paths that start at $b$. There is a commutative diagram of fibrations with contractible total spaces:
$$
\begin{CD}
 \Omega_bBG @>>> \calp_bBG @>>> BG\\
 @VV\simeq V @VV\simeq V @|  \\ 
G @>>> EG @>>> BG 
\end{CD}
$$
The homotopy equivalence on the left gives a {\em universal holonomy map}, which induces the isomorphism
\[
[X,\call BG]_b \ = \ [X, \Omega_bBG] \ \cong \ [X,G]\,.
\]
This translates the homotopy lifting action into the pointwise action of $[X,G]$ on $[X,G/T]$. 

We thus only need to identify the stabilizer from the general theory with the one claimed in Theorem 3. For this purpose, recall the homotopy equivalence of $\call BG$ with the associated bundle with respect to the conjugation action of $G$ on itself:
\[
\call BG \ \simeq \ EG \times_{(G,\conj)} G
\]
This is compatible with the projection map to $BG=EG/G$ and follows from playing with the commutative diagram above.  Applying the same discussion to $BT$, we see that the relevant map $\call\pi:\call BT\to \call BG$ can be replaced by a map
\[
EG \times_{(T,\conj)} T \  \lto \ EG \times_{(G,\conj)} G\,.
\]
Since $\pi = Bi:BT\to BG$ is induced from the inclusion $i:T\hookrightarrow G$, the map in question is simply $[\id_{EG} \times i]$, where the square bracket stands for the identifications from the group actions. Finally, since $T$ is abelian, the conjugation action is trivial and we have
\[
EG \times_{(T,\conj)} T \ \cong \ EG/T \times T \ = \ BT \times T\,.
\]

It follows that $[X, \call BT]_{\bar u} \cong [X, BT \times T]_{\bar u} \cong [X,T]$. 
Recall that $u:X\to G/T$ induces $\bar u:X\to EG/T=BT$ via $\bar u(x) = e\cdot u(x)$ where $e\in EG$ lies over $b\in BG$. (And vice versa, knowing $e$ and $\bar u$, we can recover $u$.)
For $v:X\to T$ we can now compute in $EG \times_{(G,\conj)} G \simeq \call BG$:
$$
\begin{aligned}
\call\pi \circ(\bar u \times v)(x) & \ = \ [\id_{EG} \times i]\circ(\bar u \times v)(x) 
 \ = \ [\id_{EG} \times i](e \cdot u(x), v(x)) \\
& \ = \ [e \cdot \widetilde u(x), v(x)]  
 \ = \ [e, \widetilde u(x) v(x)  \widetilde u(x)^{-1} ] \\
& \ = \ [e, u(x) v(x)  u(x)^{-1} ] 
\end{aligned}
$$
Here $\widetilde u(x)\in G$ is a momentary lift of $u(x)\in G/T$, needed to do the above manipulations. However, the element $u(x) v(x)  u(x)^{-1}\in G$ makes sense again without this choice because $T$ is abelian.  This leads to the formula in Theorem 3 for the stabilizer of the $[X,G]$-action and hence completes our proof. $\qed$

\part{Proving Theorem 2 from Theorem 3}

We conclude by indicating how the special case $(G,T) = (S^3,S^1)$ of Theorem 3 yields an alternative proof of Theorem 2.  

If we identify $S^3/S^1 \cong S^2 \cong \cpone$, $BS^1\simeq \cpinfty$ and $BS^3 \simeq \hpinfty$, then the fibration $S^3/S^1 \to BS^1 \to BS^3$ that arises in the proof of Theorem 3 becomes the familiar bundle
$$
\cpone \ \lhookrightarrow \ \cpinfty \ \overset{\pi}{\lto} \ \hpinfty
$$
where $\pi$ is induced by the inclusion $\bc\hookrightarrow\bh$.  

Since $\cpinfty$ is a $K(\bz,2)$, and $\hpinfty$ is a $5$-dimensional approximation to a $K(\bz,4)$, we have $[\bx,\cpinfty] \cong H^2(\bx)$ and $[\bx,\hpinfty] \cong H^4(\bx)$ for {\it any $4$-complex $X$}, and in this case it is also known that $\pi$ induces the cup square $Sq:H^2(\bx) \to H^4(\bx)$.  With these identifications, the map $h$ in Theorem 3 is the pull back map $h^2:\pitwo\to H^2(\bx)$, defined in the introduction.  Thus for any $u:\bx\to S^2$ (representing an element) in $\pitwo$, Theorem 3 gives an associated ``exact sequence" 
\[
\pione \ \overset{\kappa_u}{\lto}\ \pithree \ \lto\ \pitwo \ \overset{h^2}{\lto} \ H^2(\bx) \ \overset{Sq}{\lto}\ H^4(\bx)\,.
\]
Here $\kappa_u(v) = \kappa\circ(u\times v)$, where $\kappa:S^2\times S^1\to S^3$ is {\it conjugation} $(gS^1,t) \mapsto gtg^{-1}$.  The action of $\pithree$ on $\pitwo$ simply comes from the standard action of $S^3$ on $S^2=S^3/S^1$ by left translation. 

Now for $\bx$ as in Theorem 2 -- a smooth, closed, oriented and connected $4$-manifold -- the Pontryagin-Thom construction yields for any $\alphat\in\bftwo$ a dual sequence
\[
\bfthree \ \overset{\tau_\alphat}{\lto}\ \bfone \ \lto\ \bftwo \ \overset{h_2}{\lto} \ \htwo \ \overset{sq}{\lto}\ \hzero = \bz
\]
where $\tau_{\,\alphat} = \deg(\kappa)\,i_\alphat$ (here $i_\alphat$ is the framed intersection map from \secref{two})\foot{This follows from the general principle that ``suspension" in $\pi^*(X)$ corresponds to transverse ``framed intersections" in $\bff_*(\bx)$ for any $n$-manifold $\bx$.  That is, for any $p,q$ there is a commutative diagram
$$
\begin{CD}
\pi^p(X)\times\pi^q(X) @>{\sigma}>> \pi^{p+q}(X)\\
@VVV @VVV \\ 
\bff_{n-p}(X)\times\bff_{n-q}(X) @>i>>  \bff_{n-(p+q)}(X)  
\end{CD}
$$
where $\sigma$ is the {\it suspension} map, induced by any {\it degree one} map $S^p\times S^q\to S^{p+q}$, and $i$ is the {\it framed intersection} map.} 
and $sq$ is the self-intersection map.  The exactness at $\htwo$ shows that $\im(h_2) = \htwoo$, which establishes the first statement in Theorem 2, and a simple geometric argument shows that $\deg(\kappa) = 2$, establishing the last statement.  The group structure on $\bfone$, given by disjoint union, corresponds to the group structure on $\pithree$ inherited from $S^3$, as shown at the beginning of \secref{one}.

It remains to prove that the action of $\bfone$ on $\bftwo$, coming from the {\it standard action} of $S^3$ on $S^2$, corresponds to our {\it framed link translation}.   To accomplish this, view $S^3$ as the unit quaternions and $S^2$ as $S^3/S^1$, where $S^1$ is the unit circle in $\bc\subset\bh$.  Using stereographic projection, we can then visualize $S^3$ as $\br^3\cup\infty$, with $1$ at the origin, $-1$ at $\infty$, and $i,j,k$ at the tips of the standard basis vectors.  The points in $S^2$ correspond to Hopf circles in $\br^3\cup\infty$.  In particular, we will use the antipodal points $P=S^1$ (the $i$-axis) and $Q=jS^1$ (the unit circle in the $jk$-plane) in $S^2$  to analyze the standard action.

Consider a framed knot $K_e$ and a framed surface $F_\nu$ representing arbitrary elements in $\bfone$ and $\bftwo$.  Choose maps $f:\bx\to S^3$ and $g:\bx\to S^2$ with regular values $-1\in S^3$ and $Q \in S^2$ such that  $K = f^{-1}(-1)$ and $F=g^{-1}(Q)$.  We can assume that $K$ and $F$ have disjoint tubular neighborhoods $N_K$ and $N_F$ such that $f\equiv1$ off $N_K$ and $g\equiv P$ off $N_F$. Then $f\cdot g$ is equal to $g$ off $N_K$, and takes value $f(x)\cdot P$ for any $x\in N_K$.   

It follows that $Q$ is a regular value for $f\cdot g$, and that
$$
(f\cdot g)^{-1}(Q) \ = \ F \,\cup\, T
$$
where $T$ is the $2$-torus $f^{-1}(Q) \cap N_K$ (now viewing $Q\subset S^3$).  Furthermore, the framing on $T$ given by the Pontryagin construction is specified by any pushoff $f^{-1}(C)\cap N_K$, where $C$ is a Hopf circle in $S^3$ other than $P$ or $Q$.   But this is the same framed peripheral torus $T_{\tau(e)}$ for $K$ that arises in defining the action of $K_e$ on $F_\nu$ (see \figref{framing}).  Since $f\cdot g = g$ near $F$, the framing on $F$ is identical in both constructions.  This completes the proof.

\vskip-.1in

\vskip-.1in

%%%%%%%%%%%%    REFERENCES     %%%%%%%%%%%%%%%%%%%%%%%%%%%%%%%%%%%%%%%%%%%%%%%%%%%%%%%%%%%% 

\end{document}